\documentclass{article}

\usepackage[utf8]{inputenc}
\usepackage[english]{babel}

\usepackage{graphicx}
\usepackage{amssymb}
\usepackage{amsmath}
\usepackage{tikz-cd}
\usepackage{enumerate}
\usepackage{verbatim}

\usepackage{hyperref}

\newtheorem{theorem}{Theorem}[section]

\newtheorem{lemma}[theorem]{Lemma}
\newtheorem{proposition}[theorem]{Proposition}

\begin{document}
\title{Counting SET-free sets}
\author{Nate Harman}
\date{}
\maketitle

\begin{abstract}
We consider the following counting problem related to the card game SET:  \emph {How many $k$-element SET-free sets are there in an $n$-dimensional SET deck?}  Through a series of algebraic reformulations and reinterpretations, we show the answer to this question satisfies two polynomiality conditions.

\end{abstract}

\begin{section}{Introduction and formulation}

% I was told a while back that an early draft of a paper I wrote read more like I was telling someone about my research at a bar than like a research paper for publication.  Since then, I have been trying to fix this in my formal mathematical writing. This will \textbf{not} be formal mathematical writing. -- Nate

%\vspace{.7cm}

It is fairly well known in mathematical circles that the card game SET has a geometric interpretation over the field with three elements (see \cite{SET} for a non-mathematical description of the game).  Namely there is a bijection between the cards in a SET deck and $\mathbb{F}_3^4$, such that three cards form a SET if and only if the corresponding vectors $\vec{x},\vec{y},\vec{z}$ in $\mathbb{F}_3^4$ satisfy $\vec{x}+\vec{y}+\vec{z} = 0$, which geometrically means they all lie on a line.

This allows us to translate problems about the game SET into concrete mathematical questions.  For example, the question \emph{What is the maximal number of cards that can be on the table without containing a SET?} (We call such a collection a SET-free set)  translates into: \emph{What is the largest cardinality of a subset $S \subset \mathbb{F}_3^4$ such that $\nexists \ \vec{x},\vec{y},\vec{z} \in  {{S} \choose{3}}$ with $\vec{x}+\vec{y}+\vec{z} = 0$?} In this case the answer is well known to be 20.  

Such a reformulation is nice since it suggests a number of different and interesting ways we could modify or generalize the question. For example we could modify the finite field we are working over, the structure of the SETs we are looking to avoid, the dimension of the space, etc. If we just change the dimension we obtain the following famous open problem:

\medskip

\noindent \textbf{Affine Cap Problem:} What is the largest cardinality of a subset $S \subset \mathbb{F}_3^n$ such that $\nexists \ \vec{x},\vec{y},\vec{z} \in  {{S} \choose{3}}$ satisfying $\vec{x}+\vec{y}+\vec{z} = 0$?

\medskip

This is particularly nice since we can interpret it in terms of the game SET. Going from $\mathbb{F}_3^4$ to $\mathbb{F}_3^n$ corresponds to adding more attributes. In addition to looking at color, shape, number, and shading, we might also want to keep track of background color, odor, texture, political affiliation, etc. These new attributes are subject to the usual condition that there are three  possible values for each attribute, and a SET consists of three cards whose values in each attribute are either all the same or all different. We'll call the number of such attributes the \emph{dimension} of the deck. The question remains the same: \emph{How large can a SET-free set be in an $n$-dimensional SET deck?} 

This problem is well studied and notoriously difficult, and will \textbf{not} be the focus of this note.  Rather, we will relax the problem and focus on the following variation:

\medskip

\noindent \textbf{Our variation:}  In terms of $n$ and $k$, how many $k$-element subsets $S \subset \mathbb{F}_3^n$ are there with no $\vec{x},\vec{y},\vec{z} \in  {{S} \choose{3}}$ satisfying $\vec{x}+\vec{y}+\vec{z} = 0$?

\medskip

In other words, \emph {How many $k$-element SET-free sets are there in an $n$-dimensional deck?}
Or if we divide through by $\binom{3^n}{k}$, the total number of $k$-element subsets, we could ask: \emph{What is the probability of a random set of $k$-cards from an $n$-dimensional deck being SET-free?}

If we fix $n$ and ask how large $k$ can be before the answer to this problem is zero, then we recover the affine cap problem.  Instead, we will be interested in the case where $k$ is fixed, or possibly growing with $n$ but at a rate much slower than the size of an optimal solution to the affine cap problem. We will see that in this regime there is a rich algebraic structure governing the number of SET-free sets of various sizes and dimensions. %The purpose of this note will be to explore this problem a bit through a series of reformulations and applications of theorems, very little actual work is done on our end.

Let $f(n,k)$ denote the number of ordered SET-free subsets of size $k$ in an $n$ dimensional game. (Note that if we instead wanted to count unordered SET-free sets we could just divide through by $k!$.) Rather than trying to find an explicit formula for $f(n,k)$ we will focus on understanding what kind of function it is. In particular we'll see that $f(n,k)$ satisfies two polynomiality conditions.

This note is largely expository in nature and mostly consists of a series of non-technical reinterpretations of the problem, allowing us to repackage SET-free sets of various sizes and dimensions as a single algebraic object, which in turn allows us to deduce the following structural result about the function $f(n,k)$:

\begin{proposition} \label{main} Let $q = 3^n$ denote the size of the deck in an $n$-dimensional game.

\begin{enumerate}[A.]

\item If we fix $k$, then $f(n,k)$ is a polynomial of degree $k$ in $q$ with integer coefficients of alternating sign. That is, $f(n,k) = q^k - c_1(k)q^{k-1} + c_2(k)q^{k-2} - \dots \pm c_k(k)$ for some non-negative integers $c_i(k)$.

\item As we vary $k$, the coefficients $c_i(k)$ appearing in the above expression are polynomial in $k$ of degree at most $3i$.

\end{enumerate}

\end{proposition}

\end{section}

\begin{section}{Fixed $k$ and hyperplane arrangements}

Let's look at the situation in part A, where $k$ is a fixed number and we are varying the dimension $n$. We'll now make a few simple modifications which will allow us to better analyze the problem:

\begin{itemize}

\item First, let's replace $\mathbb{F}_3^n$ by $\mathbb{F}_q := \mathbb{F}_{3^n}$. They have the same additive structure, so this does not affect which subsets correspond to SETs.

\item Since we decided to look at ordered subsets, we may think of them as vectors $(x_1,x_2, \dots x_k) \in \mathbb{F}_q^k$ subject to the condition $x_i - x_j \ne 0$ for all $i,j$.  

\item Since we included an order, the condition that no three cards form a set now breaks into $\binom{k}{3}$ linear conditions $x_i + x_j + x_\ell \ne 0$ for all $i,j,\ell$.

\end{itemize}

So we have a collection of $\binom{k}{2} + \binom{k}{3}$ hyperplanes $\mathcal{H}_k = \{(x_1,\dots,x_k) \ | \ x_i - x_j = 0\}_{i,j} \cup \{ (x_1,\dots,x_k) \ | \ x_i+x_j+x_\ell = 0\}_{i,j,\ell}$ in $\mathbb{F}_q^k$, and the complement of their union counts (up to a factor of $k!$) the number of SET-free sets of $k$ cards in an $n = log_3(q)$-dimensional game of set.

%As an aside, we can make now a crude first approximation for the number of point in this complement by subtracting off $q^{k-1}$ points for each hyperplane in $\mathcal{H}_k$, ignoring intersections:

%$$ f(n,k) =  |\mathbb{F}_q^k - \mathcal{H}_k| \gtrsim q^k - \bigg(\binom{k}{2} + \binom{k}{3}\bigg)q^{k-1}$$
%In particular this shows that so long as $\binom{k}{2} + \binom{k}{3} < q$ there exist $k$-element SET-free sets. This happens at approximately $k = Cq^{1/3} = C3^{n/3}$ for some explicit constant $C$, but this isn't particularly interesting as it is what we would get with a naive probabilistic estimate (i.e. computing for what values of $k$ is the expected number of SETs at most 1).  Anyway this is off track, we were fixing $k$ and letting the dimension $n$ vary.

%\medskip

The key observation is that the actual equations defining $\mathcal{H}_k$ do not depend on $n$, and are all just defined over $\mathbb{F}_3$.  So rather than counting collections of points in spaces of different dimensions, we can now think of our problem as counting points of a single hyperplane arrangement complement (viewed as an affine scheme, if you are into that sort of thing) over different fields. This may seem like it has made the problem more complicated, but really it is application of a useful mathematical trick: By repackaging a series of related mathematical objects into a single (more complicated) one we get a more rigid theory to work with.

Hyperplane arrangements are important mathematical objects, with a fairly well developed theory. In particular the question of counting points in complements of hyperplane arrangements over finite fields has been considered by combinatorialists, and the following theorem can be found in most introductory texts  about hyperplane arrangements (see \cite{St} for example):

\begin{theorem}\label{pointcounting}
Let $\mathcal{A} \subset \mathbb{A}^k$ be a hyperplane arrangement defined over $\mathbb{F}_p$ for some prime $p$ (and hence over $\mathbb{F}_q = \mathbb{F}_{p^n}$ for all $n$).  Let $L(\mathcal{A})$ denote the lattice of flats of $A$ (viewed as subschemes of $\mathbb{A}^k$) ordered by reverse containment with minimal element $\hat{0}$ denoting the ambient space, and let $\mu$ denote the M\"obius function on $L(\mathcal{A})$.   Then:

  $$| \mathbb{F}^k_{q} - \mathcal{A}| = \chi_{\mathcal{A}}(q)$$ where $$\chi_\mathcal{A}(t) = \sum_{x\in L(\mathcal{A})} \mu(\hat{0},x)t^{dim(x)}$$ is the characteristic polynomial of the lattice of subspaces corresponding to $\mathcal{A}$ (also called the characteristic polynomial of the arrangement). Moreover, the M\"obius function satisfies $(-1)^{codim(x)} \mu(\hat{0},x) > 0$.
\end{theorem}

\noindent \textbf{Remark:} In the literature this theorem will often be stated for an arrangement $\mathcal{A}_\mathbb{Z}$ over $\mathbb{Z}$ where we are counting points of the base change $\mathcal{A}_{\mathbb{F}_q}$ to $\mathbb{F}_q$, and will require an additional assumption that we are in sufficiently large characteristic.  We'll note however that the above theorem holds (with the same proof) in any fixed characteristic, and the large characteristic assumption is only used to ensure that $L(\mathcal{A}_\mathbb{Z}) = L(\mathcal{A}_{\mathbb{F}_p})$.

\medskip

\noindent \textbf{Proof of part A of Proposition \ref{main}:} Apply of the above theorem to our arrangement $\mathcal{H}_k$ encoding the SET-free condition. $\square$

\end{section}

\begin{section}{From combinatorics to cohomology}

In order to prove part B of Proposition \ref{main}, that the coefficient $c_i(k)$ of $q^{k-i}$ appearing in $f(n,k)$ is a polynomial in $k$, we will need to again recast the problem in a different light to give an interpretation for these coefficients.

Recall that we initially defined the characteristic polynomial $\chi_\mathcal{A}(x)$ of a hyperplane arrangement $\mathcal{A}$ as a combinatorial object, namely as the characteristic polynomial of a certain lattice associated to $\mathcal{A}$. This is a nice definition and works well over for hyperplane arrangements over arbitrary fields, however for a complex hyperplane arrangement $\mathcal{A} \subset \mathbb{C}^k$ there is a celebrated theorem of Orlik and Terao which expresses $\chi_\mathcal{A}(x)$ in terms of the topology of the complement:

\begin{theorem}{\label{OTthm}}\textbf{(\cite{OT} Theorem 3.68)}
$$\chi_\mathcal{A}(x) = \sum (-1)^i x^{k-i} \text{dim} \ H^i(\mathbb{C}^k - \mathcal{A},\mathbb{Q})$$
where $H^i(\mathbb{C}^k - \mathcal{A},\mathbb{Q})$ is the singular cohomology of the hyperplane complement. 
\end{theorem}
So if we were working over $\mathbb{C}$ these coefficients we want to understand would just be the Betti numbers of the corresponding hyperplane complements.

\medskip

However we aren't working over $\mathbb{C}$, we working over finite fields, and singular cohomology is not a useful notion. Lucky for us, modern (well, like 50 years old at this point) mathematics has already solved this whole ``singular cohomology doesn't make sense for varieties over finite fields" issue with the invention of \'etale cohomology. Indeed, if $\mathcal{A} \subset \mathbb{A}^k$ is a hyperplane arrangement defined over $\mathbb{F}_p$ and $\ell \ne p$ is a prime number we have the following version of the above theorem due to Kim \cite{Kim}:
$$\chi_\mathcal{A}(x) = \sum (-1)^i x^{k-i} \text{dim} \ H_{\acute{e}t}^i(\mathbb{A}^k - \mathcal{A},\mathbb{Q}_\ell).$$ 

As an aside, we'll note here that evaluating at $x= q := p^n$ as in Theorem \ref{pointcounting} to count points can be viewed as an application of the Grothendieck-Lefschetz trace formula, but we won't go into that in any more detail.
\medskip

So in this light, we see that these individual coefficients are really counting something topological about these hyperplane arrangements.  In particular we see that part B of Proposition \ref{main} is equivalent to the following lemma, the proof of which we will leave to the next section. 

\begin{lemma}\label{mainthm}
$\dim H_{\acute{e}t}^i(\mathbb{A}^k - \mathcal{H}_k,\mathbb{Q}_\ell)$ is a polynomial in $k$ of degree at most $3i$.

\end{lemma}

\end{section}

\begin{section}{$FI$-modules and polynomiality in $k$}

The purpose of this section will be to prove Lemma \ref{mainthm}.  To do so we will once again apply our favorite trick of repackaging a series of related objects into a single object with more structure.  In this case we will replace the series of hyperplane arrangements $\mathcal{H}_k$ by a single object $\mathcal{H}$, something called an $FI$-CHA.  As a result the series of cohomology spaces $H_{\acute{e}t}^i(\mathbb{A}^k - \mathcal{H}_k,\mathbb{Q}_\ell)$ (for fixed $i$, $k$ varying) will automatically be given the structure of a finitely generated free $FI$-module.

Let's unravel this a bit:  $FI$ stands for the category of $F$inite sets with $I$njections, an $FI$-module (over a field $\mathbb{F}$, for simplicity) is a covariant functor from $FI$ to vector spaces over $\mathbb{F}$.  What this amounts to is a representation $V_k$ of the symmetric group $S_k$ for each $k$ along with maps $V_j \to V_k$ for every injection $[j] \to [k]$ satisfying some natural compatibility constraints.  

The point is that (under some reasonable technical conditions) the sequences of representations $V_k$ that can show up in this setting are far from arbitrary. We have the following theorem (with all technical conditions surpressed):

\begin{theorem}\label{FIthm}\textbf{(Summarizing \cite{CEF1} and \cite{CEFN})} \

\begin{enumerate}[(1)]
\item Under a finite generation condition, the sequence of dimensions $\text{dim} \ V_k$ is eventually polynomial in $k$.

\item In characteristic zero and under the same condition as (1), the sequence of representations $V_k$ exhibits the phenomenon known as representation stability. -- We won't use this, but we'll mention it since it's really cool.

\item Under the same condition as (1) and an additional freeness condition, the dimensions $\text{dim} \ V_k$ are actually polynomial in $k$ (as opposed to eventually polynomial).

\end{enumerate}
\end{theorem}

\medskip

An $FI$-CHA or $FI$-Complement of Hyperplane Arrangement is a collection of hyperplane arrangements $\mathcal{A}_k \subset \mathbb{A}^k$ (for our purposes defined over $\mathbb{F}_p$ for some $p$) for each $k\in \mathbb{N}$ such that:

\begin{enumerate}[(a)]
\item $\mathcal{A}_k$ contains the hyperplane defined by $x_1-x_2 = 0$.

\item $\mathcal{A}_k$ is preserved by the ``permute the coordinates" action of $S_k$ on $\mathbb{A}^k$

\item The inverse image of a hyperplane in $\mathcal{A}_j$ under the ``forget the last $k-j$ coordinates" map $\mathbb{A}^k \to \mathbb{A}^j$ is a hyperplane in $\mathcal{A}_k$.

\item There is a finite collection of hyperplanes such that every hyperplane in $\mathcal{A}_k$ for all $k$ can be obtained by repeated applications of rules (2) and (3).

\end{enumerate}

The point of this definition is that the hyperplane arrangement complements along with various ``forget and reorder some of the coordinates" maps define a well behaved contravariant functor from $FI$ to a certain category of smooth schemes over $\mathbb{F}_p$.  Taking \'etale cohomology (which is again contravariant) gives the following theorem:

\begin{theorem}\label{FICHA} \textbf{(\cite{CEF2} Corollary 3.2 and \cite{Gadish} Theorem 1.4)}
Let $\mathcal{A}$ be an $FI$-CHA over $\mathbb{F}_p$, and let $\ell \ne p$ be a prime number. For each $i$, the spaces $H_{\acute{e}t}^i(\mathbb{A}^k - \mathcal{A}_k,\mathbb{Q}_\ell)$ form a finitely generated free (in the sense of Theorem \ref{FIthm}) $FI$-module.

\end{theorem}

We are now ready to prove Lemma \ref{mainthm} (and hence part B of Prop. \ref{main}):

\medskip

\noindent \textbf{Proof of Lemma \ref{mainthm}:} Combining Theorems \ref{FIthm} and \ref{FICHA}, in order to see that the dimensions are polynomial it suffices to verify that the SET arrangements form a $FI$-CHA. So for that we just need to check the four conditions above:

\begin{itemize}
\item The $S_k$ action in (b) just comes from the fact that we passed from sets to ordered sets in section 2, but a SET-free set is SET-free no matter how we order it.

\item Condition (a), along with the $S_k$ action encodes the statement that sets do not have repeated elements.

\item The finite list in condition (d) is just $x_1+x_2+x_3=0$ (along with $x_1-x_2 =0$, which we just addressed) which up to symmetry encodes the SET-free condition.

\item Condition (c) just means that a subset of a SET-free set is SET-free.

\end{itemize}

We still need to prove the degree bound for the polynomials, that is, that $c_i(k)$ is a polynomial of degree at most $3i$ in $k$. In principle we could do this by going through the above calculations of $FI$-modules more carefully to obtain this bound directly.

 Instead we will again apply our favorite repackaging trick and recall that if we take the cohomology spaces $H_{\acute{e}t}^i(\mathbb{A}^k - \mathcal{H}_k,\mathbb{Q}_\ell)$ for various $i$ together they form a graded ring. Moreover, for hyperplane arrangements this ring is known to be generated in degree  1 (this is another fairly well known fact from \cite{OT}).

In particular this implies that it is enough to check the degree bound in the case when $i=1$. By looking at the lattice of flats associated to $\mathcal{H}_k$ it is easy to calculate that:
$$\dim H_{\acute{e}t}^1(\mathbb{A}^k - \mathcal{H}_k,\mathbb{Q}_\ell) = \binom{k}{2} + \binom{k}{3}$$
which is indeed degree $3$ in $k$. $\square$

\end{section}

\begin{section}{Generalizations and problems}

We chose to focus on the SET example mostly because it is fun, but pretty much everything we did here works more or less the same if do one or more of the following:

\begin{itemize}

\item Work in fixed characteristic $p \ne 3$.

\item Modify the linear constraint we are avoiding, say $\vec{x} + \vec{y} = \vec{z}$ or $2\vec{x} + 2\vec{y} = \vec{z}$ instead of $\vec{x} + \vec{y} + \vec{z} = 0$.

\item Avoid multiple (but finitely many) linear constraints.

\item Avoid higher codimension linear constraints, such as $\vec{x},\vec{y},\vec{z},\vec{w}$ forming a four term arithmetic progression. (Here we need a generalization of $FI$-CHAs found in \cite{Gadish})
\end{itemize}

In all these cases we will again get that the number of $k$ element subsets of $\mathbb{F}_{p}^n$ avoiding these constraints will be polynomial in $q=p^n$ for fixed $k$, and that the coefficients of these polynomials will vary polynomially in $k$.  Again we'll stress that we don't expect this to give any useful insights into the corresponding extremal problems.

\medskip

\noindent \textbf{Problem 1: (Mixed characteristics)} In many of the cases we care about, the equations are all defined over $\mathbb{Z}$.  For a fixed hyperplane arrangement with equations defined over $\mathbb{Z}$ the combinatorics of the arrangement over $\mathbb{C}$ are the same as the combinatorics over the arrangement defined over $\mathbb{F}_p$ for all but finitely many primes $p$.

This breaks down in general for an $FI$-CHA, which consists of infinitely many hyperplane arrangements.  However we conjecture that the additional ``bad" primes we get by looking at $\mathcal{A}_k$ for larger and larger $k$ will only affect the cohomology in larger and larger degrees.  More precisely, we conjecture that for fixed $i$ the polynomial $d_i(k) =  \text{dim } H_{\acute{e}t}^i(\mathbb{A}^k - \mathcal{A}_k, \mathbb{Q}_\ell)$ is independent of characteristic away from finitely many primes.

\medskip

\noindent \textbf{Problem 2: (Computations)} We might wonder what these polynomials actually are. Computing by hand, the first few terms of $f(n,k)$ are:

$$ f(n,k) = q^k - \bigg(\binom{k}{2} + \binom{k}{3}\bigg)q^{k-1} + \bigg(2\binom{k}{3}+ 15\binom{k}{4} + 25\binom{k}{5} + 10 \binom{k}{6} \bigg) q^{k-2}  + O(q^{k-3})$$
It should certainly be possible to generate more of them using a computer, and perhaps say more structurally about the polynomials appearing.  More ambitiously, we could ask what are the stable multiplicities of irreducible symmetric group representations coming from the $FI$-module structure (the dimensions of which correspond to these formulas)? 

\medskip

\noindent \textbf{Problem 3: (Analyzing convergence)} For fixed $k$ these formulas are just polynomial in $q$ so there are no convergence issues. However if we think of $k$ as growing with $q$ these formulas involve more and more terms and we have to deal with convergence issues.  If we normalize by dividing by $\binom{q}{k}$ and ask about the probability of a random $k$-element subset being SET-free, then a reasonable first guess might be that if we let $k$ grow proportional to $q^c$ for $c < \frac{1}{3}$ then this probability tends to $1$. Is this true? What can be said about the rate of convergence? What happens as $c$ goes to $1/3$? Are there interesting scaling limits?

\end{section}

\end{document}